\documentclass[11pt]{article}
\usepackage{amssymb}

\def\<{\langle}
\def\>{\rangle}

\def\ol{\overline}
\def\wh{\widehat}
\def\a{\alpha}
\def\b{\beta}

\def\l{\lambda}

\def\Q{\mathbb{Q}}
\def\Z{\mathbb{Z}}
\def\R{\mathbb{R}}
\def\C{\mathbb{C}}
\def\L{\Lambda}
\def\G{\Gamma}
\def\D{\Delta}
\def\nsg{\triangleleft}
\def\mm#1#2#3#4{\left(\begin{array}{cc}#1 &#2\\ #3 &#4\end{array}\right)}

\newtheorem{thm}{Theorem}[section]
\newtheorem{lem}[thm]{Lemma}
\newtheorem{cor}[thm]{Corollary}
\newtheorem{conj}{Conjecture}

\title{Generalised triangle groups of type $(3,5,2)$}
\author{James Howie\\
Department of Mathematics and\\ Maxwell Institute for Mathematical Sciences\\
Heriot-Watt University\\
Edinburgh EH14 4AS\\
UK}

\begin{document}

\maketitle

\begin{abstract}
If $G$ is a group with a presentation of the form
$\< x,y|x^3=y^5=W(x,y)^2=1\>$, then either $G$ is
virtually soluble or $G$ contains a free subgroup of rank $2$.
This provides additional evidence in favour of a conjecture
of Rosenberger.
\end{abstract}

\section{Introduction}

A {\em generalised triangle group} is a group $G$ with a presentation of
the form $$\< x,y|x^p=y^q=W(x,y)^r=1\>$$ where $p,q,r\ge 2$ are integers
and $W(x,y)$ is a word of the form $$W(x,y)=x^{\a(1)}y^{\b(1)}\cdots x^{\a(k)}y^{\b(k)},$$
with $0<\a(i)<p$ and $0<\b(i)<q$.  We say that $G$ is of {\em type} $(p,q,r)$.
%The parameter $k$ is called the {\em length}.
Without loss of generality,
we assume that $p\le q$.

A conjecture of Rosenberger \cite{Ros} asserts that a Tits alternative holds for
generalised triangle groups:

\begin{conj}[Rosenberger]\label{RC}
Let $G$ be a generalised triangle group.  Then
either $G$ is soluble-by-finite
or $G$ contains a non-abelian free subgroup.
\end{conj}

In a  recent article \cite{Howie2}, we proved
the Rosenberger Conjecture in the case $(p,q,r)=(3,3,2)$.
In the present note we prove it
in the case $(p,q,r)=(3,5,2)$.
 In conjunction with previously
known results \cite{FLR,BMS,Howie,LR,BK,BK2,BKB,BKB1,BKB2,HW,HW2,Will2}
(see for example the survey \cite{FRR} or the Introduction
to \cite{Howie2} for details), this reduces the conjecture
to the cases $(p,q,r)=(2,q,2)$ for $q\in\{3,4,5\}$.

\section{Preliminary results}

Suppose that $X,Y\in SL(2,\C)$ are matrices, and $W=W(X,Y)$
is a word in $X,Y$.  Then the trace of $W$ can be calculated
as the value of a $3$-variable polynomial, where the variables
are the traces of $X$, $Y$ and $XY$ \cite{Horowitz}.  We can use
this to find and analyse {\em essential representations}
from $G$ to $PSL(2,\C)$.  (A representation of $G$ is {\em essential}
if the images of $x,y,W(x,y)$ have orders $p,q,r$ respectively.)

We can force the images $x,y$ to have orders $3,5$ in $PSL(2,\C)$
by mapping them to matrices $X,Y\in SL(2,\C)$ of trace $2\cos(\pi/3)=1$
and $2\cos(\pi/5)=(1+\sqrt{5})/2$ respectively.  Then the trace of
$W(X,Y)\in SL(2,\C)$ is given by a one-variable polynomial $\tau_W(\l)$
of degree $k$,
where $\l$ denotes the trace of $XY$.  We obtain an essential representation by
choosing $\l$ to be a root of $\tau_W$ (which forces the image of $W$ to have
order $2$ in $PSL(2,\C)$.

\begin{lem}\label{roots}
$G$ contains a nonabelian free subgroup, unless the roots of $\tau_W(\l)$
all belong to\\ $\{0,1,(1+\sqrt{5})/2,(-1+\sqrt{5})/2\}$.
\end{lem}

\noindent{\em Proof}.
The image of an essential representation of $G$ is generated by two elements
of orders $3$ and $5$ respectively, and contains an element of order $2$.
With the exception of the finite group $A_5$, any such subgroup of
$PSL(2,\C)$ contains a nonabelian free subgroup.  The result follows, unless all
essential representations $\rho:G\to PSL(2,\C)$ have image isomorphic to $A_5$.

Now let $\l$ be a root of $\tau_W(\l)$ corresponding to an essential representation
$\rho:G\to A_5$, where $\rho(x),\rho(y)$ are represented by
matrices $X,Y$ of traces $1,(1+\sqrt{5})/2$ respectively.  Then
$XY$ and $XY^{-1}$ are matrices representing nontrivial elements
of $A_5$, which therefore have orders in $\{2,3,5\}$.  Thus the traces of
$XY$ and $XY^{-1}$ belong to $\{0,\pm 1,(\pm 1\pm\sqrt{5})/2\}$.  Moreover,
these traces  also satisfy the trace equation
$$tr(XY)+tr(XY^{-1})=tr(X)tr(Y)=\frac{1+\sqrt{5}}{2}.$$
From this, it follows that $\l=tr(XY)\in\{0,1,(\pm 1+\sqrt5)/2\}$, as claimed.

\begin{lem}\label{l1}
Let $p:\ol{K}\to K$ be a regular covering of connected
$2$-complexes with $K$ finite, with covering transformation
group abelian of torsion-free rank at least $2$.  Let $F$ be a field.  If
$$H_2(\ol{K},F)=0\ne H_1(\ol{K},F),$$
then
$$\mathrm{dim}_F~H_1(\ol{K},F)=\infty.$$
\end{lem}

\noindent{\em Proof}.
Let $\{a,b\}$ be a basis for a free abelian subgroup $A$ of the
group of covering
transformations of $p:\ol{K}\to K$, and let $\alpha$ be a
cellular $1$-cycle of $\ol{K}$ over $F$ that represents a non-zero
element of $H_1(\ol{K},F)$.  If the $F[a]$-submodule of $H_1(\ol{K},F)$
generated by $\alpha$ is free, then $H_1(\ol{K},F)$ is infinite-dimensional
over $F$,
as claimed.  So we may assume that there is a cellular $2$-chain $\beta$
of $\ol{K}$ with $d(\beta)=f(a)\alpha$ for some non-zero
polynomial $f(a)\in F[a]$.

For similar reasons, we may also assume that $d(\gamma)=g(b)\alpha$
for some cellular $2$-chain $\gamma$ of $\ol{K}$
and some non-zero polynomial $g(b)\in F[b]$.

Now $f(a)\gamma-g(b)\beta\in H_2(\ol{K},F)=0$.  In other words
$f(a)\gamma=g(b)\beta$ in the group $C_2(\ol{K},F)$ of cellular $2$-chains
of $\ol{K}$, which is  a free
module over the unique factorisation domain $FA\cong F[a^{\pm 1},b^{\pm 1}]$.
Since $f(a),g(b)$ are coprime in $F[a^{\pm 1},b^{\pm 1}]$, it follows
that there is a $2$-chain $\delta$ with $f(a)\delta=\beta$ and $g(b)\delta=\gamma$.
Hence $f(a)(d(\delta)-\alpha)=d(\beta)-f(a)\alpha=0$, in the group $C_1(\ol{K},F)$ of
cellular $1$-chains of $\ol{K}$.  But $C_1(\ol{K},F)$ is also
a free module over the domain $F[a^{\pm 1},b^{\pm 1}]$,
and $f(a)\ne 0$, so $d(\delta)=\alpha$, contradicting the hypothesis that
$\alpha$ represents a non-zero element of $H_1(\ol{K},F)$.

This contradiction completes the proof.

\begin{lem}\label{l2}
 Let $E$ be the set of midpoints of edges of a
regular icosahedron $\mathcal{I}\subset\R^3$ centred at the origin, and let
$M=\Z E$ its $\Z$-span in $\R^3$.
Let $V=\{1,a,b,c\}\subset Isom^+(\mathcal{I})\subset SO(3)$ be the Klein $4$-group,
and let $C=\{1,c\}\subset V$. Then, regarding $M$ as a $\Z V$-module via
the action of $V$ by isometries of $\mathcal{I}$, we have the following.
\begin{enumerate}
\item $M\cong\Z^6$ as an abelian group.
\item $H_0(C,M)=\Z\otimes_{\Z C} M\cong\Z_2^4\oplus\Z^2$.
\item The induced action of $V/C$ on $H_0(C,M)/(\mathrm{torsion})$ is mutliplication by
$-1$.
\end{enumerate}
\end{lem}

\noindent{\em Proof}.
If $e$ is the midpoint of the edge joining two vertices
$u,v$ of $\mathcal{I}$, then $e=(u+v)/2$.  Thus $E$ is contained in the $\Q$-span
$W$ of the set of vertices of $\mathcal{I}$.  Since the vertices occur in $6$ antipodal
pairs, the $\Q$-span $\Q M$ of $E$ has dimension at most $6$ over $\Q$.

On the other hand, for any vertex $v$, $\sqrt{5}\cdot v$ is the sum of
the $5$ vertices adjacent to $v$ in $\mathcal{I}$.  Thus $\sqrt{5}\cdot v\in W$.
It also follows that $\sqrt{5}\cdot e\in M$ for any $e\in E$:
specifically, $(\sqrt{5}-2)\cdot e$ is the sum of the midpoints
of the eight edges of $\mathcal{I}$
that share a vertex with the edge containing $e$.  If $e_1,e_2,e_3\in E$
are chosen to be linearly independent over $\R$ -- and hence over $\Q[\sqrt{5}]$ --
then $e_1,e_2,e_3,\sqrt{5}\cdot e_1, \sqrt{5}\cdot e_2, \sqrt{5}\cdot e_3\in M$
are linearly independent over $\Q$.  Thus $\Q M=\Q\otimes_\Z M$ has dimension
exactly $6$ over $\Q$.  Since $M\subset\Q M$ is torsion-free and finitely
generated, it follows that $M\cong\Z^6$, as claimed.

If, in the above, we choose $e_1,e_2,e_3$
to lie on the axes of the rotations
$a,b,c\in V$ respectively, 
then we obtain a decomposition
$$\Q M=\Q[\sqrt{5}]e_1 \oplus \Q[\sqrt{5}]e_2 \oplus \Q[\sqrt{5}]e_3$$
of $\Q M$ as a $\Q[\sqrt{5}]$-vector space, with respect to which $a,b,c$ act as the diagonal matrices
$diag(1,-1,-1)$, $diag(-1,1,-1)$ and $diag(-1,-1,1)$
respectively.  Let $$M_+:=M\cap\Q[\sqrt{5}]e_3~~~\mathrm{and}~~~ 
M_-:=M\cap (\Q[\sqrt{5}]e_1 \oplus \Q[\sqrt{5}]e_2).$$  Then $M_-\cap M_+=\{0\}$,
while $e_1,e_2,\sqrt{5}e_1,\sqrt{5}e_2\in M_-$ and $e_3,\sqrt{5}e_3\in M_+$,
so $M_-$, $M_+$ are free abelian of ranks $4$ and $2$ respectively.

Moreover, $M/M_-$ is naturally embedded in the vector space $\Q M/\Q M_-$, so
is also free abelian -- necessarily of rank $2$.  Note that $M_-$ is closed under the
action of $V$ on $M$.  Under the induced action on $M/M_-$, each of $a,b$ acts as
the antipodal map, multiplication by $-1$, and $c$ acts as the identity.

Hence $(1-c)M=2M_-$, so $$H_0(C,M)=M/(1-c)M=M/2M_-\cong \Z_2^4\oplus\Z^2,$$
as claimed.

Finally, the quotient of $H_0(C,M)$ by its torsion subgroup is naturally
isomorphic to $M/M_-$, and the induced action of $V/C$ on this quotient is
via the antipodal map.

\begin{lem}\label{l3}
Let $G=\<x,y|x^3=y^5=W(x,y)^2=1\>$ and suppose that
$(\lambda-\a)^2$ divides the trace polynomial $\tau_W(\lambda)$
of $W$, for some $\a\in\{0,1,(1+\sqrt{5})/2,(-1+\sqrt{5})/2\}$.
Let $\rho:G\to%\<x,y|x^3=y^5=(xy)^2=1\>\cong 
A_5$
be the natural epimorphism corresponding to the root $\a$ of $\tau_W(\l)$.
Let $C\subset A_5$ be a subgroup
of order $2$ and $V\subset A_5$ its centraliser of order $4$.
Then $G$ has subgroups $N_1\nsg N_2\nsg\rho^{-1}(V)$
such that
\begin{enumerate}
\item $\rho(N_2)=\{1\}$;
\item $\rho^{-1}(C)/N_2\cong\Z^2$;
\item $\rho^{-1}(V)/N_2\cong\Z^2\rtimes_{(-1)}\Z_2$;
\item $N_2/N_1$ is a non-zero vector space over $\Z_2$.
\end{enumerate} 
\end{lem}

\noindent{\em Proof}.
Let $\L=\mathbb{C}[\lambda]/\<(\lambda-\a)^2\>$, and choose matrices
$$X=\mm{e^{i\pi/3}}{0}{1}{e^{-i\pi/3}},~Y=\mm{e^{i\pi/5}}{\l-\a-2\cos(8i\pi/15)}{0}{e^{-i\pi/5}}\in SL_2(\L)$$ so that $$tr(X)=1,~~tr(Y)=\frac{1+\sqrt{5}}{2}~~
\mathrm{and}~~tr(XY)=\lambda-\a.$$  Then $X,Y$ determine a representation
$\wh\rho:G\to PSL_2(\L)$, since $tr(W(X,Y))=\tau_W(\l)=0$ in $\L$.  If $\phi:PSL_2(\L)\to PSL_2(\mathbb{C})$ is the
natural epimorphism obtained by setting $\lambda=\a$, then
the image of $\rho=\phi\circ\wh\rho$ is isomorphic to $A_5$.
Let $K$ denote the kernel of $\rho$ and let $L$ denote
the kernel of $\wh\rho$.

Clearly $G/K \cong A_5$.  Now $K/L\cong\wh\rho(K)$ is the normal closure
of $(xy)^2.L$, so it is isomorphic to the subgroup of $PSL(2,\L)$
generated by
$$(XY)^2=-I+(\lambda-\a)(XY)$$
together with its conjugates by elements of $\wh\rho(G)$.
Let $Z=\phi(XY)\in A_5\subset SU(2)$ denote the matrix obtained from $XY$
by subsituting $\lambda=\a$. Note that $tr(Z)=0$, in other
words, $Z\in sl_2(\mathbb{C})$.  Since $(\lambda-\a)^2=0$ in $\L$, we also have
$$(XY)^2=-I+(\lambda-\a)Z.$$
For similar reasons,
for any $M\in\wh\rho(G)$ we have
$$M(XY)^2M^{-1}=-I+\phi(M)Z\phi(M)^{-1}.$$
Moreover, since $(\lambda-\a)^2=0$ in $\L$ we have, for any $A,B\in sl_2(\mathbb{C})$,
$$(-I+(\lambda-\a) A)(-I+\lambda B)=I+(\lambda-\a) (A+B).$$
Thus $K/L\cong\rho(K)$ is isomorphic to the additive subgroup
of $sl_2(\mathbb{C})$ generated by $MZM^{-1}$ for all $M\in \widehat{A_5}\subset SU(2)$.
There are precisely $30$ such conjugates of $Z$; geometrically they correspond to the
midpoints of the edges of a regular icosahedron centred at the origin in $\R^3$,
where we identify $SU(2)$ with the $3$-sphere of unit-norm quaternions,
and $\R^3$ with the space of purely imaginary quaternions.  As an abelian
group, therefore, $K/L\cong\rho(K)\cong\mathbb{Z}^6$ by Lemma \ref{l2}.

Now $K/L$ is also an $A_5$-module.  Its structure as an $A_5$-module does not need to concern
us, but Lemma \ref{l2} gives us some information about its structure as a $C$-module
and as a $V$-module.  This in turn gives information on the structure
of $\D:=(\rho)^{-1}(C)$.

Specifically, $H_0(C,K/L)=H_0(\D/K,K/L)\cong\Z_2^4\oplus\Z^2$.  It follows from the
$5$-term exact sequence
$$H_2(\D/L)\to H_2(\D/K)\to H_0(\D/K,K/L)\to H_1(\D/L)\to H_1(\D/K)\to 0$$
and the fact that $\D/K\cong\Z_2$
that $H_1(\D/L)$ has torsion-free rank $2$, and that the torsion subgroup
of $H_1(\D/L)$ is a non-zero vector space over $\Z_2$.  Hence we can define $N_1=[\D,\D].L$ and $N_2\supset N_1$
such that $N_2/N_1$ is the torsion-subgroup of $\D/N_1=H_1(\D/L)$.
That $N_1\nsg \rho^{-1}(V)$ follows from the fact that $[\D,\D]$ and $L$ are
both normal in $\rho^{-1}(V)$.  That $N_2~\nsg~\rho^{-1}(V)$ follows
from the fact that $N_2/N_1$ is characteristic in $\D/N_1$.

Finally, since $V/C$ acts on $\Z^2\cong \D/N_2$ by the antipodal map, it follows
that $\rho^{-1}(V)/N_2\cong\Z^2\rtimes_{(-1)}\Z_2$, as required.

\section{Main results}

\begin{thm}\label{mult}
Let $G=\<x,y|x^3=y^5=W(x,y)^2=1\>$.  If the
trace polynomial $\tau_W(\lambda)$ of $W$ has a multiple root, then
$G$ contains a nonabelian free subgroup.
\end{thm}

\noindent{\em Proof}.
We may assume that the root $\a$ is one of $0$,
$1$, $(\pm 1+\sqrt{5})/2$, for otherwise the result is
immediate from Lemma \ref{roots}.  Let $\rho:G\to A_5$
be the essential representation corresponding to $\a$,
let $c=\rho(W)\in A_5$, $C=\{1,c\}\subset A_5$ the subgroup generated
by $c$, and $V=\{1,a,b,c\}\subset A_5$ its centraliser in $A_5$.

Let $N_1\nsg N_2\nsg \rho^{-1}(V)<G$ be the subgroups promised by Lemma \ref{l3}.
Let $\G=\rho^{-1}(C)<\rho^{-1}(V)$ be the unique index $2$ subgroup such that $N_2\subset\G$
and $\G/N_2\cong\Z^2$.  Then $\G$ has index $30$ in $G$ and contains
no conjugate of $x$ or of $y$.

Applying the Reidemeister-Scheier process to the presentation
of $G$ in the statement of the Theorem, we obtain a presentation
of $\G$ of the form
$$\G=\< k_1,\cdots,k_{31}|r_1,\dots,r_{30},s_1^2,s_2^2\>,$$
where $r_1,\dots,r_{10}$ are rewrites of conjugates
of $x^3$; $r_{11},\dots,r_{16}$ are rewrites
of conjugates of $y^5$; and $r_{17},\dots,r_{30}$ and $s_1^2=W^2,s_2^2=\hat{a}W^2\hat{a}^{-1}$
are rewrites of conjugates of $W^2$, with $\rho(\hat{a})=a$ and
so $s_1=W,s_2=\hat{a}W\hat{a}^{-1}\in\G$.

Let $K$ be the $2$-complex model of this presentation,
$F=\Z_2$, and $p:\ol{K}\to K$ the regular cover correspdonding
to the normal subgroup $N_2\nsg\G$.  Let $L\subset K$
be the subcomplex obtained by omitting the $2$-cells
corresponding to the relators $s_1^2,s_2^2$, and let
$\ol{L}:=p^{-1}(L)\subset\ol{K}$.

Now, since $\G/N_2$ is torsion-free, and since $s_1^2=1=s_2^2$
in $\G$, $s_1,s_2\in N_2$.  Hence each lift of each $2$-cell
$s_i^2$ ($i=1,2$) to $\ol{K}$ is bounded by the square of some
path in $\ol{K}^{(1)}$.  As a consequence, the $2$-cells in
$\ol{K}\setminus\ol{L}$ represent elements
of $H_2(\ol{K},F)$, and it follows that the inclusion-induced
map $H_1(\ol{L},F)\to H_1(\ol{K},F)$ is an isomorphism.

Since $N_2/N_1$ is a nonzero
$F$-vector space, we have $$H_1(\ol{L},F)\cong H_1(\ol{K},F)=H_1(N_2,F)\ne 0.$$

If $H_2(\ol{L},F)=0$,
then
by Lemma \ref{l1} it follows that $dim_F~H_1(N_2,F)=\infty$.
On the other hand, if $H_2(\ol{L},F)\ne 0$ then $H_2(\ol{L},F)$
contains a free $F(\G/N_2)$-module of rank $>0=\chi(L)$, since $F(\G/N_2)$
is an integral domain.  In this case $H_1(\ol{L},F)$ contains a non-zero
free $F(\G/N_2)$-submodule, by \cite[Proposition 2.1 and Theorem 2.2]{Howie}.  Again we deduce that
$dim_F~H_1(N_2,F)=\infty$.

Thus the Bieri-Strebel invariant $\Sigma$ of the $F(\G/N_2)$-module
$N_2/N_1$ is a proper
subset of $S^1$ \cite[Theorem 2.4]{BS}.  But by Lemma \ref{l2} (3) it follows that $\Sigma$
is invariant under the antipodal map: $\Sigma=-\Sigma$.
Hence $\Sigma\cup -\Sigma\ne S^1$, and it follows \cite[Theorem 4.1]{BS} that $\G$
contains a nonabelian free subgroup,
as claimed.

\begin{cor}[Main Theorem]
Let $G$ be a generalised triangle group of type $(3,5,2)$.
Then either $G$ is virtually soluble or $G$ contains a nonabelian
free subgroup.
\end{cor}

\noindent{\em Proof}.
By Theorem \ref{mult} and Lemma \ref{roots} the result follows
unless $\tau_W(\l)$ has only simple roots in the set
$\{0,1,(1+\sqrt{5})/2,(-1+\sqrt{5})/2\}$, in which case the
degree $k$ of $\tau_W(\l)$ is at most equal to $4$.

But the Rosenberger Conjecture is known for $k\le 4$
\cite{LR}.

\end{document}